\documentclass[a4paper,12pt]{amsart}
\usepackage{amssymb}
\usepackage{ifthen}
\usepackage{graphicx}
\usepackage{color}
\nonstopmode
\newtheorem{thm}{Theorem}
\newtheorem{cor}{Corollary}
\newtheorem{lem}{Lemma}

\newtheorem{claim}{Claim}

\newtheorem{conj}{Conjecture}
\newtheorem{prob}{Problem}

\theoremstyle{definition}
\newtheorem{defn}{Definition}
\newtheorem{example}{Example}

\newenvironment{rem}{%
\bigskip
\noindent \textsl{{\sl Remark. }}}{\bigskip}
\newenvironment{rems}{%
\bigskip
\noindent \textsl{{\sl Remarks. }}}{\bigskip}

\newenvironment{pf}[1][]{%
 \vskip 1mm
 \noindent
 \ifthenelse{\equal{#1}{}}%
  {{\slshape Proof. }}%
  {{\slshape #1.} }%
 }%
{\qed\medskip}

\newcounter{alphabet}

\newcommand{\ID}{{\mathbb D}}

\newcommand{\D}{{\mathbb D}}
\def\be{\begin{equation}}
\def\ee{\end{equation}}

\newcommand{\bee}{\begin{enumerate}}
\newcommand{\eee}{\end{enumerate}}

\newcommand{\blem}{\begin{lem}}
\newcommand{\elem}{\end{lem}}
\newcommand{\bthm}{\begin{thm}}
\newcommand{\ethm}{\end{thm}}
\newcommand{\bcor}{\begin{cor}}
\newcommand{\ecor}{\end{cor}}
\newcommand{\beg}{\begin{example}}
\newcommand{\eeg}{\end{example}}
\newcommand{\begs}{\begin{examples}}
\newcommand{\eegs}{\end{examples}}
\newcommand{\bdefe}{\begin{defn}}
\newcommand{\edefe}{\end{defn}}
\newcommand{\bprob}{\begin{prob}}
\newcommand{\eprob}{\end{prob}}
\newcommand{\bques}{\begin{ques}}
\newcommand{\eques}{\end{ques}}
\newcommand{\bei}{\begin{itemize}}
\newcommand{\eei}{\end{itemize}}

\newcommand{\bde}{\begin{deter}}
\newcommand{\ede}{\end{deter}}
\newcommand{\bca}{\begin{case}}
\newcommand{\eca}{\end{case}}
\newcommand{\bcl}{\begin{claim}}
\newcommand{\ecl}{\end{claim}}
\newcommand{\bcon}{\begin{conj}}
\newcommand{\econ}{\end{conj}}
\newcommand{\bcons}{\begin{conjs}}
\newcommand{\econs}{\end{conjs}}
\newcommand{\bprop}{\begin{propo}}
\newcommand{\eprop}{\end{propo}}
\newcommand{\br}{\begin{rem}}
\newcommand{\er}{\end{rem}}
\newcommand{\brs}{\begin{rems}}
\newcommand{\ers}{\end{rems}}
\newcommand{\bo}{\begin{obser}}
\newcommand{\eo}{\end{obser}}
\newcommand{\bos}{\begin{obsers}}
\newcommand{\eos}{\end{obsers}}
\newcommand{\bpf}{\begin{pf}}
\newcommand{\epf}{\end{pf}}
\newcommand{\ba}{\begin{array}}
\newcommand{\ea}{\end{array}}
\newcommand{\beq}{\begin{eqnarray}}
\newcommand{\beqq}{\begin{eqnarray*}}
\newcommand{\eeq}{\end{eqnarray}}
\newcommand{\eeqq}{\end{eqnarray*}}

\newcounter{minutes}
\divide\time by 60
\newcounter{hours}
\multiply\time by 60 \addtocounter{minutes}{-\time}
\begin{document}
\title[An elementary counterexample to a coefficient conjecture  ]
{An elementary counterexample to a coefficient conjecture }
%[Geometric studies on the family $\mathcal{G}(\alpha)$]
%{Geometric studies on the family $\mathcal{G}(\alpha)$}

%\title[Generalized Zalcman conjecture for a family  of univalent functions]
%{Generalized Zalcman conjecture for  a family  of univalent
%functions}

%=========================================================================
\thanks{%$^\dagger$
File:~\jobname .tex,
          printed: \number\day-\number\month-\number\year,
          \thehours.\ifnum\theminutes<10{0}\fi\theminutes}
%=========================================================================

\author{Liulan Li, Saminathan  Ponnusamy and Karl-Joachim Wirths
%$^\dagger $
%${}^{~\mathbf{*}}$
}
\address{L. Li, College of Mathematics and Statistics
 (Hunan Provincial Key Laboratory of Intelligent Information Processing and Application),\\
 Hengyang Normal University,
 Hengyang,  Hunan 421002, People's Republic of China.}
\email{lanlimail2012@sina.cn}

\address{S. Ponnusamy,  Department of Mathematics,
Indian Institute of Technology Madras, Chennai-600 036, India.}
\address{Department of Mathematics, Petrozavodsk State University, ul., Lenina 33, 185910 Petrozavodsk, Russia}
\email{samy@iitm.ac.in}

\address{K.-J. Wirths, Institute of Algebra and Analysis,
Technical University Braunschweig, 38106 Braunschweig, Germany.}
 \email{k-j.wirths@tu-bs.de}


%\\
%$%{}^{\mathbf{*}}

%^\dagger$

%\thanks{ }
%\maketitle

%\newfont{\Bbb}{msbm10 scaled\magstephalf}
\begin{abstract}
In this article, we consider the family of  functions $f$ meromorphic in the unit disk $\ID=\{z :\,|z| < 1\}$ with a pole at
the point $z=p$, a Taylor expansion
\[f(z)= z+\sum_{k=2}^{\infty} a_kz^k, \quad |z|<p,
\]
and satisfying the condition
\[\left |\left(\frac{z}{f(z)}\right)-z\left(\frac{z}{f(z)}\right)'-1\right |<\lambda,\, \forall z\in\ID,
\]
for some $\lambda$, $0<\lambda < 1$. We denote this class by $\mathcal{U}_m(\lambda)$ and we shall prove a representation theorem for the functions in this class.
As consequences, we get a simple proof for the estimates of $|a_2|$ and obtain inequalities for the initial coefficients of the Laurent series of
$f\in \mathcal{U}_m(\lambda)$ at its pole. In \cite{PW2} it had been conjectured that for $f\in \mathcal{U}_m(\lambda)$ the inequalities
\[|a_n|\,\leq\,\frac{1}{p^{n-1}}\sum_{k=0}^{n-1}(\lambda p^2)^k, \quad n\geq 2
\]
are valid. We provide a counterexample to this conjecture for the case $n=3$.
\end{abstract}

\maketitle
\pagestyle{myheadings}
\markboth{L. Li, S. Ponnusamy and K.-J. Wirths}{A counterexample}

\section{Introduction}

Let $\lambda \in (0,1)$, $p\in (0,1),$ and $\mathcal{U}_m(\lambda)$
be the family of functions meromorphic in the open unit disk
$\ID=\{z :\,|z| < 1\}$ with a pole at the point $z=p$,  a Taylor
expansion \be\label{f1} f(z)= z+\sum_{k=2}^{\infty} a_kz^k, \quad
|z|<p, \ee and satisfying the condition \be\label{fb}
|U_f(z)|<\lambda, ~\mbox{ for every}~ z\in\ID, \ee where
\[U_f(z)\,:=\,\left(\frac{z}{f(z)}\right)-z\left(\frac{z}{f(z)}\right)'-1.
\]
It has been proved in \cite{A} that the functions in
$\mathcal{U}_m(\lambda)$ are univalent in $\ID$. For a different
proof in the related case of  holomorphic functions $f$ in $\ID$,
please refer to \cite{NO}.
%The function $f\in
%\mathcal{U}_m(\lambda)$ has at most one pole in $\ID$, which is
%necessarily a simple pole, see \cite{A}.

For $f\in \mathcal{U}_m(\lambda)$ with the expansion \eqref{f1}, in
\cite{PW2} it has been proved by an application of Banach's fixed
point theorem that the sharp inequality \be\label{fa}
|a_2|\,\leq\,\frac{1+\lambda p^2}{p} \ee is valid, and it has been
conjectured there that the inequalities \be\label{f2}
|a_n|\,\leq\,\frac{1}{p^{n-1}}\sum_{k=0}^{n-1}(\lambda p^2)^k
~\mbox{ for  $n\geq 2$}, \ee hold. In this paper, we shall prove a
representation formula for the functions in $\mathcal{U}_m(\lambda)$
which take the pole $z=p$ into account. By using the above
representation, we will give a simple proof of the inequality
\eqref{fa}, and obtain estimates for the initial coefficients of the
Laurent series of the functions in $\mathcal{U}_m(\lambda)$ at the
pole. Further we shall show that the inequality \eqref{f2} is not
true for $n=3$.

%\bigskip
%\noindent
As the title indicates, the main aim of this paper is to disprove the conjecture \eqref{f2}. We want to emphasize at this place that some of the other results can be found in \cite{BP1, BP2,BP3}. As a byproduct of main concern, we include some easy consequences of our discussion
in the form of corollaries.

\section{Representation and estimates of coefficients}

Let ${\mathcal B}$ denote the set of the functions $\omega$ which are analytic
in $\D$ and satisfy $|\omega(z)|\leq 1$ for all $z\in\D$.

\bthm\label{representation} {\rm (Compare \cite{BP2})} The function $f$ belongs to the class $\mathcal{U}_m(\lambda)$ and has a simple pole at the point $z=p$ if and only if
there exists a function $\omega\in {\mathcal B}$ such that
\be\label{fc}
f(z)\,=\,\frac{z}{1\,-\,\frac{z}{p}\,+\,\lambda z \int_z^p\omega(t)\,dt}.
\ee
\ethm

\bpf Necessity. Since $f\in \mathcal{U}_m(\lambda)$ has a simple pole at the point $z=p$, we use the normalization at the origin to
get the representation
\[\frac{z}{f(z)}\,=\,1\,-\,\frac{z}{p}\,+\,r(z),
\]
where the function $r$ is holomorphic in $\ID$ and satisfies the condition $r(p)=0.$ A calculation yields
\[U_f(z)\,=\,r(z)\,-\,zr'(z).
\]
On the other hand we see that $U_f(z)$ has a zero of multiplicity at least two at the origin. This fact together with the inequality
\eqref{fb} yields the equation
\[U_f(z)\,=\,\lambda z^2 \omega(z), \quad z\in \ID,
\]
where $\omega\in {\mathcal B}$. Now, it is easy to find the solution for the initial value problem
\[r(z)\,-\,zr'(z)\,=\,\lambda z^2 \omega(z),\quad r(p)=0.
\]
Obviously it is given by
\[r(z)\,=\,\lambda z \int_z^p\omega(t)\,dt.
\]

Sufficiency is obvious. \epf
\medskip

By Theorem \ref{representation}, we have the following corollaries.

\bcor {\rm (Compare \cite{BP1})} Let $f\in \mathcal{U}_m(\lambda)$ be as in the form of
\eqref{fc} with a simple pole at the point $z=p$ and a Taylor
expansion as in \eqref{f1}. Then the region of variability of the
coefficient $a_2$ is the disk
\[\left\{\frac{1}{p}\,-\,\lambda p u\,:\, u\in
\overline{\ID}\right\},
\]
and
\[\frac{1}{p}\,-\,\lambda p\,\leq\,|a_2|\,\leq\, \frac{1}{p}\,+\,\lambda p,
\]
where the upper bound is just \eqref{fa}.
\ecor

\bpf  By Theorem \ref{representation} we get
\[a_2\,=\,\frac{1}{p}\,-\,\lambda\int_0^p\omega(t)\,dt.
\]
Hence, the conclusions hold.\epf

\bcor\label{b-1} {\rm (Compare \cite{BP3})} Let $f\in \mathcal{U}_m(\lambda)$ be as in the form
of \eqref{fc} with a simple pole at the point $z=p$, where the
Laurent series of $f$ at $z=p$ is
\be\label{laurent series}
f(z)\,=\,\sum_{k= -1}^{\infty}b_k(z-p)^k.
\ee
Then the region of variability of the residue $b_{-1}$ is the disk
\[\left\{\frac{-p^2}{1\,+\,\lambda p^2 u}\,:\, u\in \overline{\ID}\right\},\]
and
\[\frac{p^2}{1\,+\,\lambda p^2}\,\leq\,|b_{-1}|\,\leq\, \frac{p^2}{1\,-\,\lambda p^2},\]
where the upper bound is sharp and the sharpness is attained when
$\omega (z)\equiv\,-1, z\in \ID$.
\ecor

\bpf A little calculation using the above theorem yields that
\be\label{fA} b_{-1}\,=\,\frac{-p^2}{1\,+\,\lambda p^2
\omega(p)}.\ee Then the conclusion is obvious. \epf

\blem\label{the function of p} Let $p\in \big (\frac{\sqrt{17}-1}{4}, 1\big )$
and $\phi(p)=\frac{2p^2+p-2}{p^3}$. Then $\phi(p)\in(0,1).$\elem

\bpf Set $t=1/p$ and $\psi(t)=-2t^3+t^2+2t$. Then
$t\in \big (1,\frac{\sqrt{17}+1}{4}\big )$ and $\phi(p)=\psi(t)$. A calculation
yields that
$$\psi'(t)=-6t^2+2t+2=-6\left(t-\frac{\sqrt{13}+1}{6}\right)\left(t+\frac{\sqrt{13}-1}{6}\right)<0$$
for $t\in \big (1,\frac{\sqrt{17}+1}{4}\big )$. Thus, $\psi'(t)$ is
monotonically decreasing for $t\in \big (1,\frac{\sqrt{17}+1}{4}\big )$ and
$$0=\psi \Big (\frac{\sqrt{17}+1}{4}\Big )
<\psi(t)<\psi(1)=1.$$ The proof is completed. \epf

\bcor\label{b0} Let $f\in \mathcal{U}_m(\lambda)$ be as in the form
of \eqref{fc} with a simple pole at the point $z=p$, where the
Laurent series of $f$ at $z=p$ is as in the form of \eqref{laurent
series}. Then we have the following conclusions.
\begin{enumerate}

\item[{\rm (i)}] If $p\in \big (0,\frac{\sqrt{17}-1}{4} \big]$, then
\[|b_0|\leq \frac{p}{(1\,-\,\lambda p^2)^2}.\]

\item[{\rm (ii)}] If $p\in \big (\frac{\sqrt{17}-1}{4}, 1 \big)$ and
$\lambda\in \big [\frac{2p^2+p-2}{p^3},1\big )$,
 then
\[|b_0|\leq \frac{p}{(1\,-\,\lambda p^2)^2}.\]

\item[{\rm (iii)}] If $p\in \big(\frac{\sqrt{17}-1}{4}, 1\big )$ and
$\lambda\in \big(0,\frac{2p^2+p-2}{p^3}\big),$
 then
\[|b_0|\leq \frac{p}{2(1\,-\, p^2)}\frac{\lambda p^3-2p^2+2}{1-\lambda p(\lambda p^3-2p^2+2)}.\]

\end{enumerate}

The equalities in {\rm (i)} and {\rm (ii)} occur for $\omega(z)\equiv -1,\,z\in
\ID,$ while the equality in {\rm (iii)} is attained for the case
\[\omega(z)\,=\,-\,\frac{a\,+\,z}{1\,+\,az},\quad z\in \ID,
\]
where $a\in (-p,1)$ is chosen such that
\[  \frac{a\,+\,p}{1\,+\,ap}\,=\,\frac{2\,-\,2p^2\,+\,\lambda p^3}{p}.
\]\ecor

\bpf  A calculation using the above theorem and Corollary \ref{b-1}
yields that \be\label{fB} b_0\,=\,\frac{-2p\,+\lambda
p^4\omega'(p)}{2(1\,+\,\lambda p^2\omega(p))^2}. \ee

To get a sharp upper bound for $|b_0|$, we use the well known
inequality
\[|\omega'(z)|\,\leq\,\frac{1\,-\,|\omega(z)|^2}{1\,-\,|z|^2},\quad z\in
\ID.
\]
Hence, \eqref{fB} together with the last inequality implies
\[|b_0|\,\leq\,\frac{p}{2(1-p^2)}\,\frac{2\,-\,2p^2\,+\,\lambda p^3\,-\lambda p^3|\omega(p)|^2}{(1\,-\,\lambda p^2|\omega(p)|)^2}.
\]
For further calculations we use the following abbreviations
\[ A\,=\,2\,-\,2p^2\,+\,\lambda p^3,~~\,B\,=\,\lambda p^3,\,C\,=\,\lambda p^2,\]
and
\[ x\,=\,|\omega(p)|\in [0,1],\quad D(x)\,=\,\frac{A\,-\,B x^2}{(1\,-\,Cx)^2}.\]
Then
\[D'(x)\,=\,2\,\frac{-B x\,+\,CA}{(1\,-\,C x)^3}=\frac{2\lambda p^3}{(1\,-\,C
x)^3}\left(\frac{\,2\,-\,2p^2\,+\,\lambda p^3}{p}-x\right).
\]
\medskip

\noindent{\bf (i)} If $p\in \big (0,\frac{\sqrt{17}-1}{4} \big ]$, then
$$2p^2+p-2=2\left(p-\frac{\sqrt{17}-1}{4}\right)\left(p+\frac{\sqrt{17}+1}{4}\right)\leq0,
$$
and thus,
$$\lambda p^3-(2p^2+p-2)>0, \; \ ~\mbox{ i.e. }~\;\ \frac{\,2\,-\,2p^2\,+\,\lambda p^3}{p}>1.
$$
Therefore, we get a maximum of $D(x)$ at $x=1$ and (i) holds.
\medskip

If $p\in \big (\frac{\sqrt{17}-1}{4}, 1 \big)$, then Lemma \ref{the function of
p} implies that $\frac{2p^2+p-2}{p^3}\in(0,1).$\medskip

\noindent{\bf (ii)} If $p\in \big(\frac{\sqrt{17}-1}{4}, 1 \big )$ and
$\lambda\in \big[\frac{2p^2+p-2}{p^3},1 \big),$ then
$$\frac{\,2\,-\,2p^2\,+\,\lambda p^3}{p}>1.
$$
Therefore, we get a maximum of $D(x)$ at $x=1$ and (ii) holds.
\medskip

\noindent{\bf (iii)} If $p\in \big(\frac{\sqrt{17}-1}{4}, 1 \big)$ and
$\lambda\in \big(0,\frac{2p^2+p-2}{p^3} \big),$ then
$$0<\frac{\,2\,-\,2p^2\,+\,\lambda p^3}{p}<1.
$$
Therefore, we get a maximum of $D(x)$ at
$x=\frac{\,2\,-\,2p^2\,+\,\lambda p^3}{p}$ and (iii) holds.\epf

\br In \cite{BP3}, Bhowmik and Parveen gave Conjecture 1 that the
Laurent coefficients $b_n$ of $f\in \mathcal{U}_m(\lambda)$ at $z=p$
satisfy
$$|b_n|\leq\frac{\lambda^n p^{n+1}}{(1-\lambda p^2)^{n+2}},\,
n\geq0,$$ and the equality above holds for the function
$f(z)=pz/(p-z)(1-\lambda pz).$ For the case $n=0$, Corollary
\ref{b0} shows that this conjecture holds if $p\in \big
(0,\frac{\sqrt{17}-1}{4} \big]$, or $p\in \big
(\frac{\sqrt{17}-1}{4}, 1 \big)$ and $\lambda\in \big
[\frac{2p^2+p-2}{p^3},1\big )$ while for $p\in
\big(\frac{\sqrt{17}-1}{4}, 1\big )$ and $\lambda\in
\big(0,\frac{2p^2+p-2}{p^3}\big),$ Corollary \ref{b0} shows that
this conjecture does not hold. \er

\section{A counterexample}

Our counterexample to the above mentioned conjecture from \cite{PW2} is as follows and it generalizes a counterexample from \cite{LPW} to
a conjecture from \cite{OPW}. To prove the example, we need the following lemma.

\blem\label{g(p)}
Let $g(p)=4\ln(1+p)\,-\,3p$. Then there exists a unique $p_0\approx 0.7336$ such that
\begin{enumerate}
\item[{\rm (i)}]  $g(p)>0$ for $p\in(0,p_0)$;

\item[{\rm (ii)}]  $g(p_0)=0$;

\item[{\rm (iii)}] $g(p)<0$ for $p\in(p_0,1)$.
\end{enumerate}
\elem
\bpf Since $g'(p)=\frac{1-3p}{1+p}$, it is easy to see that the function $g(p)$ is monotonically increasing for $p\in \big [0,\frac{1}{3}\big]$ and monotonically decreasing
for $p\in \big (\frac{1}{3},1\big ]$. Since $g(0)=0$, $g(p)>0$ for $p\in \big(0,\frac{1}{3}\big]$. Since $g(1)=4\ln2-3<0$, $g(1/3)>g(0)=0$ and $g(p)$ is monotonically decreasing for
$p\in \big (\frac{1}{3},1 \big]$, there exists a unique $p_0$ such that $g(p_0)=0$ and $g(p)<0$ for $p\in(p_0, 1]$. A calculation yields $p_0\approx 0.7336$.
\epf

\beg\label{a counterexample} Let $p_0\approx 0.7336$ be the unique solution of the equation $4\ln(1+p)\,-\,3p=0$ for $p\in(0,1)$. For any $p\in(p_0,1)$ and
\[\lambda \in\left(0,\frac{\frac{4}{p}\log(1+p)-3}{-4p\log(1+p)+2p^2}\right),
\]
there exists an $f\in \mathcal{U}_m(\lambda)$ such that the third coefficient $a_3(f)$ of the Taylor expansion of $f$ in $|z|<p$ satisfying
$$|a_3(f)|\,>\frac{1}{p^2}\left(1+\lambda p^2\,+\,\lambda^2p^4\right).
$$
\eeg

\bpf Let $a \in [0,1]$ and
\be\label{f3}
f_a(z)\,=\,\frac{z}{1-\frac{z}{p}-\lambda z  \int_z^p\frac{a+t}{1+at}\,dt}.\
\ee
It is clear from Theorem \ref{representation} that $f_a \in \mathcal{U}_m(\lambda)$. Set
\[ v_p(a)\,:=\,\int_0^p\frac{a+t}{1+at}\,dt.
\]
Then from \eqref{f3}, the denominator of $f_a$ is as follows.
\[1-\frac{z}{p}-\lambda z\int_z^p\frac{a+t}{1+at}\,dt\,=\,1-z\left(\frac{1}{p}+\lambda v_p(a)\right)\,+\,\lambda a z^2\,+\,o(z^2).
\]
This implies the expansion
\[f_a(z)\,=\, z\,+\,z^2\left(\frac{1}{p}+\lambda v_p(a)\right)\,+\,z^3\left(\left(\frac{1}{p}+\lambda v_p(a)\right)^2\,-\,\lambda a\right)\,+\,o(z^3).
\]
Therefore, the third Taylor coefficient of $f_a$ is given by
\[a_3\,=\,\frac{1}{p^2}\left(1\,+\,\lambda(2pv_p(a)-p^2a)\,+\,\lambda^2p^2v_p^2(a)\right).
\]
Since
\[ v_p'(a)\,=\,\int_0^p\frac{1-t^2}{(1+at)^2}\,dt\,>\,0,
\]
the function $v_p$ is monotonically increasing and $v_p(a)<v_p(1)=p$ for $a\in [0,1)$. Hence,
\be\label{f4}
a_3\,>\,\frac{1}{p^2}\left(1+\lambda p^2\,+\,\lambda^2p^4\right)
\ee
if and only if
\be\label{f5}
0\,<\,\lambda\,<\,\frac{\frac{2v_p(a)}{p}\,-\,a\,-1}{p^2\,-\,v_p^2(a)},\quad a\in(0,1).
\ee

Let us consider now the function
\[ w_p(a)\,:=\,\frac{2v_p(a)}{p}\,-a.
\]
Then $w_p(0)=p$, $w_p(1)=1$, and
\[ w_p'(a)\,=\,\frac{2v_p'(a)}{p}\,-1
\]
so that
\[ w_p'(0)\,=\,\frac{2v_p'(0)}{p}\,-1=\frac{2}{p}\left (p-\frac{p^3}{3}\right )-1 =  1-\frac{2}{3}p^2
\]
and
\[ w_p'(1)\,=\,\frac{2v_p'(1)}{p}\,-1=\frac{2}{p}\int_0^p\frac{1-t}{1+t}\,dt -1= \frac{4}{p}\ln(1+p)\,-\,3.
\]
As
\[w_p''(a)\,=\,-\frac{4}{p}\int_0^p\frac{t(1-t^2)}{(1+at)^3}\,dt\,<\,0,\]
the function $w_p'$ is monotonically decreasing from
$w_p'(0)=1-\frac{2}{3}p^2>0$ to
$$w_p'(1)=\frac{4}{p}\ln(1+p)\,-\,3.$$
Since $p\in(p_0,1)$, Lemma \ref{g(p)} implies that $w_p'(1)<0$. The
above facts yield that $w_p'(a)=0$ for $a\in (0,1)$ has the unique
solution $a_p$. Then $w_p$ is increasing from $w_p(0)=p<1$ to
$w_p(a_p)$, and decreasing from $w_p(a_p)$ to $w_p(1)=1.$ Since
\[\max\{w_p(a): a\in [0,1]\}\,=\,w_p(a_p)\,>\,1,\]
we see that for any $p\in(p_0,1)$, there exists a unique $a_0\in
(0,a_p)$ such that $$w_p(a_0)-1\,=\,0$$ and $w_p(a)-1>0$ for all
$a\in(a_0,1)$.

Therefore, we conclude that for any $p\in (p_0,1)$ and all $a\in
(a_0,1)$ we have
\[0\,<\,\frac{w_p(a)-1}{p^2-v_p^2(a)}.\]
To get a concrete interval of values of $\lambda$ wherein \eqref{f4} is satisfied, we observe that
\[L_p(a)\,:=\,\frac{w_p(a)-1}{p^2-v_p^2(a)}\]
which vanishes at $a=a_0$, is continuous on $[a_0,1)$ and has the
limit
\[\lim_{a\to 1}\frac{w_p'(1)}{-2v_p'(1)v_p(1)}\,=\,\frac{\frac{4}{p}\log(1+p)-3}{-4p\log(1+p)+2p^2}.\]
Therefore, at least for $p\in (p_0,1)$ and
\[\lambda \in\left(0,\frac{\frac{4}{p}\log(1+p)-3}{-4p\log(1+p)+2p^2}\right),\]
there exist functions $f_a$ of the above form such that \eqref{f4} is
fulfilled.\epf

\subsection*{Acknowledgements}
The work of the first author is supported by NSF of Hunan (No.
2020JJ6038), the Scientific Research Fund of Hunan Provincial
Education Department (20A070), the Scientific Research Fund of
Hengyang Normal University(No. KYZX21002), the Science and
Technology Plan Project of Hunan Province (No. 2016TP1020), and the
Application-Oriented Characterized Disciplines, Double First-Class
University Project of Hunan Province (Xiangjiaotong [2018]469).

\subsection*{Statements and Declarations}
%\subsection*{Conflict of Interests}
The authors declare that there is no conflict of interests regarding the publication of this paper.

%\subsection*{Data Availability Statement}
The authors declare that this research is purely theoretical and does not associate with any datas.

\end{document}